\documentclass[1p,times,sort]{elsarticle}
\usepackage[utf8]{inputenc}
\usepackage{amsmath,amsfonts,amssymb,amsthm}
\usepackage[all]{xy}
\usepackage{graphicx,subfig,aliascnt}
\usepackage[colorlinks=true]{hyperref}

\journal{arXiv}

\theoremstyle{plain}
\newtheorem{theorem}{Theorem}[section]
\newtheorem{proposition}[theorem]{Proposition}
\newtheorem{corollary}[theorem]{Corollary}
\newtheorem{lemma}[theorem]{Lemma}
\theoremstyle{definition}
\newtheorem{definition}[theorem]{Definition}

\newtheorem{remark}[theorem]{Remark}
\newtheorem{notation}[theorem]{Notation}

\theoremstyle{plain}
\newtheorem*{theorem*}{Theorem}
\newtheorem*{proposition*}{Proposition}
\newtheorem*{corollary*}{Corollary}
\newtheorem*{lemma*}{Lemma}
\theoremstyle{definition}
\newtheorem*{definition*}{Definition}
\newtheorem*{remark*}{Remark}
\newtheorem*{notation*}{Notation}

\newcommand{\PP}{\mathbb{P}}
\newcommand{\CC}{\mathbb{C}}

\newcommand{\ZZ}{\mathbb{Z}}

\newcommand{\QQbar}{\overline{\mathbb{Q}}}

\DeclareMathOperator{\res}{res}

\DeclareMathOperator{\spec}{Spec}

\begin{document}
\begin{frontmatter}
\title{The Casas-Alvero conjecture for three recycled roots in degree 20}

\author[dm]{C\'esar Massri\corref{correspondencia}\fnref{financiado}}
\ead{cmassri@caece.edu.ar}
\address[dm]{Department of Mathematics, CAECE, Buenos Aires, Argentina}
\cortext[correspondencia]{Address for correspondence: Department of Mathematics, CAECE, Buenos Aires, Argentina.
Postal address: Av. de Mayo 866. Phone number: 54-11-5217-7878.}
\fntext[financiado]{The author was fully supported by CONICET, IMAS, Buenos Aires, Argentina}

\begin{abstract}
The Casas-Alvero conjecture says that a 
degree $n$ complex univariate polynomial 
sharing a root with each of its derivative
must have only one root.
In this article we give three results. The first one, 
is that the number of possible counterexamples in normal form
of degree $p^r+p^s$ or $p^r+2p^s$
is finite  ($p$ prime, $r,s$ positive integers).
The second result is that a possible counterexample in normal form
of degree $p^r+1$ has algebraic
coefficients and the final result is that in degree $20$ there are no
counterexamples with three recycled roots.
\end{abstract}

\end{frontmatter}

\section*{Introduction}

E. Casas-Alvero conjectured in 
relation to his work on plane curves \cite{casas,casas2}, that
a degree $n$ polynomial $f$ sharing a root with each of its derivative
must have only one root.
Specifically, let $f\in\CC[x]$ be 
a univariate polynomial such that  
\[
\res(f,f^{i})=0,\quad 1\le i\le n-1,
\]
where $\res(-,-)$ is the resultant.
The Casas-Alvero conjecture says that $f$ 
has a multiplicity $n$ root.
This conjecture was checked with a computer in low 
degrees \cite{diaz2005conjecture}
and proved, by using techniques from number theory for
degrees $p^e,2p^e,3p^e,4p^e,5p^e$, where 
$e$ is a natural number and 
$p$ runs through a list of infinitely many prime numbers,
see \cite{verhoek2009some,chellali,draisma,infinity} for a precise
statement. Also, in \cite{yaku,yaku2,castryck} several constrains on a possible
counterexample are given. The conjecture is known to be 
false in positive characteristic and the first open case is $n=20$.

In this article, we present new results regarding the Casas-Alvero conjecture.
We based our analysis in standard arguments appearing in the literature
and on the properties of Abel-Gontcharoff polynomials.
Specifically, we applied arguments on positive characteristics and brute-force computations to deduce our first two results.
Our third result were deduced from properties of Abel-Gontcharoff polynomials
where we proved that there are no counterexamples in degree 20 with three recycled
roots.

In the first four sections we 
collect some results from the available bibliography (published or not)
and give some remarks. We do not take credit on these results except maybe
on Proposition \ref{igualdad} that we do not found it in the literature on
interpolation theory.
In section \ref{first} we prove that there are a finite number of possible
counterexamples (written
in normal form) of degrees $p^s(p^r+1)$ or $p^s(p^r+2)$, where $p$ is a prime
and $r,s$ non-negative integers.
In section
\ref{second} we prove that the possible counterexamples in degree $p^r+1$ are in $\QQbar[x]$
and finally, in section \ref{third} we prove that there are no counterexamples with
three recycled roots in degree $20$.

\section{Polynomials with one root}

\begin{definition}
Let $f$ be a degree $n$ complex polynomial and let $h$ be some
complex number. We denote $f^h$ to the polynomial given by 
\[
f^h(x):=f(x+h),\quad x\in\CC.
\]
Notice that if $a$ is a root of $f$, then $a-h$ is  a root of $f^h$.
\end{definition}

\begin{lemma}\label{equivs}
Let $f$ be a degree $n$ complex monic polynomial and let $F$ be its companion matrix. 
Then, the following are equivalent
\begin{enumerate}
\item $f$ is equal to $b(x-a)^n$ for some complex numbers $a,b$, $b\neq 0$.
\item $f$ has exactly $1$ root.
\item $\res(f,f^h)\ne0$ for all non-zero complex number $h$.
\item $\det(f^h(F))\ne0$ for all non-zero complex number $h$.
\item $\det(f^h(F))$ as a polynomial in $h$ has exactly $1$ root.
\item $\det(f^h(F))$ as a polynomial in $h$ is equal to $h^{m}$ for some $m>0$.
\item $\det(f^h(F))$ as a polynomial in $h$ is equal to $h^{n^2}$.
\item $\res(f,f^h)=h^{n^2}$.
\end{enumerate}
\end{lemma}
\begin{proof}
Let us prove $2$ implies $3$. If $f$ has exactly one root $a$, then $f^h$ also has exactly one root 
equal to $a-h$. Then, 
$\res(f,f^h)$ is zero if and only if $h=0$. Now let us prove $8$ implies $1$.
From $8$ we deduce that $\res(f,f^h)$ is zero if and only if $h=0$ and this implies $1$.
\end{proof}

\begin{lemma}
Let $f$ be a degree $n$ complex polynomial with roots $\{\lambda_1,\ldots,\lambda_n\}$
and let $F$ be its companion matrix. 
Then, $\det(f^h(F))$ as a polynomial in $h$ is equal to either of the following expressions
\[
\sum_{k=0}^{n^2}\left(\sum_{k_1+\dots+k_n=k}\frac{f^{(k_1)}(\lambda_1)}{k_1!}\dots \frac{f^{(k_n)}(\lambda_n)}{k_n!} \right)h^k=
h^n\prod_{i<j}\left(h^2 - (\lambda_i-\lambda_j)^2\right).
\]
Both expressions are symmetric polynomials in $\{\lambda_1,\ldots,\lambda_n\}$.
\end{lemma}
\begin{proof}
Given that $f^h$ is a polynomial and that $F$ is similar to a lower triangular matrix, 
we can compute the expression $\det(f^h(F))$
from the eigenvalues of $F$ which are the roots of $f$. Then, 
\[
\det(f^h(F))=f^h(\lambda_1)\dots f^h(\lambda_n) = f(\lambda_1+h)\dots f(\lambda_n+h).
\]
Recall that 
if $\lambda_1,\dots,\lambda_n$ are the roots of $f$, then the roots of 
$f(\lambda_k+h)$ are $\lambda_1-\lambda_k,\dots,\lambda_n-\lambda_k$. Hence, 
we can factorize the right hand side of the previous equation as
\begin{align*}
f(\lambda_1+h)\dots f(\lambda_n+h) &= \left( (h-\lambda_1+\lambda_1)\dots(h-\lambda_1+\lambda_n)\right)
\dots
\left((h-\lambda_n+\lambda_1)\dots(h-\lambda_n+\lambda_n)\right)\\
&= h^n\prod_{i<j}\left(h^2 - (\lambda_i-\lambda_j)^2\right)
\end{align*}
where, in the last equality, we simplified $h-\lambda_k+\lambda_k$ and 
collected the terms $(h-\lambda_i+\lambda_j)(h-\lambda_j+\lambda_i)$. 
Now, expanding in $h$ the left
hand side of the previous equation, we get
\[
f(\lambda_1+h)\dots f(\lambda_n+h) = \prod_{i=1}^n\left(\sum_{k_i=0}^n \frac{f^{(k_i)}(\lambda_i)}{k_i!}h^{k_i}\right)=
\sum_{k=0}^{n^2}\left(\sum_{k_1+\dots+k_n=k}\frac{f^{(k_1)}(\lambda_1)}{k_1!}\dots \frac{f^{(k_n)}(\lambda_n)}{k_n!} \right)h^k.
\]
\end{proof}

\begin{remark}
Let $e_k$ be the $k$-th elementary symmetric polynomial in $n(n-1)/2$ variables and degree $k$
evaluated at $\{(\lambda_i-\lambda_j)^2\}_{i<j}$.
Then,
\[
h^n\prod_{i<j}\left(h^2 - (\lambda_i-\lambda_j)^2\right)=
h^{n^2}+e_1h^{n^2-2}+e_2h^{n^2-4}+\dots+e_kh^{n^2-2k}+\dots+e_{\frac{n(n-1)}{2}}h^{n}.
\]
In particular, if $r$ is such that $n^2-r$ is odd,
\[
\sum_{k_1+\dots+k_n=r}\frac{f^{(k_1)}(\lambda_1)}{k_1!}\dots \frac{f^{(k_n)}(\lambda_n)}{k_n!} =0.
\]
Clearly, if some $k_i=0$, then $f^{(k_i)}(\lambda_i)=f(\lambda_i)=0$. Hence, 
the first non-trivial term is $f'(\lambda_1)\dots f'(\lambda_n)$, the coefficient of $h^n$.
%
\end{remark}

\section{CA-polynomials}

\begin{definition}
A degree $n$ complex univariate polynomial $f$
is called a \emph{CA-polynomial} if the
resultant between $f$ and its $i$-th derivative vanishes,
\[
\res(f,f^{(i)})=0,\quad 1\le i\le n-1.
\]
\end{definition}

\begin{notation}
Let $\epsilon_i$ be the $i$-th elementary symmetric polynomial in $n$ variables and degree $i$.
The following notation will be used throughout the article.
Let us denote $\epsilon_i^k$ to the polynomial 
$\epsilon_i(\lambda_1-\lambda_k,\dots,\lambda_n-\lambda_k)$
in $\mathbb{Z}[\lambda_1,\dots,\lambda_n]$,
\[
\epsilon_i^k:=
\epsilon_i(\lambda_1-\lambda_k,\dots,\lambda_n-\lambda_k).
\]
\end{notation}

\begin{proposition}\label{ca-equiv}
Let $f$ be a CA-polynomial with roots $\lambda_1,\dots,\lambda_n$.
The following are equivalent
\begin{enumerate}
\item $f$ is a CA-polynomial.
\item The roots of $f$ satisfies the polynomial equations,
\[
f^{(i)}(x_1)\dots f^{(i)}(x_n)=0,\quad i=1,\dots,n-1.
\]
\item For all $i=1,\dots,n-1$,
\[
\epsilon_{n-i}^1\dots \epsilon_{n-i}^n= 0.
\]
\item There exists a function 
$\varphi:\{1,\dots,n-1\}\to\{1,\dots,n\}$ such that
\[
\epsilon_{1}^{\varphi(1)}=\dots
=\epsilon_{n-1}^{\varphi(n-1)}=0.
\]
\end{enumerate}
\end{proposition}
\begin{proof}
Let us prove the equivalence between $1$ and $4$.
Let $f$ be a CA-polynomial. Given that $\res(f,f^{(i)})=0$, 
there exists some root $\lambda_k$ such that $f^{(i)}(\lambda_k)=0$. Hence, 
writing 
\[
f(\lambda_k+h)=
f(\lambda_k)+hf^{(1)}(\lambda_k) + \frac{h^2}{2!}f^{(2)}(\lambda_k)+\dots+\frac{h^n}{n!}f^{(n)}(\lambda_k),
\]
we get that the coefficient of the monomial $h^i$ is zero. But, given that the roots
of $f(h+\lambda_k)$ are $\{\lambda_1-\lambda_k,\dots,\lambda_n-\lambda_k\}$, 
it follows that the coefficient of $h^i$ is equal to $\epsilon_{n-i}^{k}$. Hence, 
\[
\epsilon_{n-i}(\lambda_1-\lambda_k,\dots,\lambda_n-\lambda_k)=0.
\]

Now, let $f$ be a degree $n$ complex polynomial and assume that there
exists a function 
$\varphi:\{1,\dots,n-1\}\to\{1,\dots,n\}$ such that
$\epsilon_i(\lambda_1-\lambda_{\varphi(i)},\dots,\lambda_n-\lambda_{\varphi(i)})=0$
for $1\le i\le n-1$. Then, the polynomial $f(\lambda_{\varphi(i)}+h)$ 
do not have the monomial $h^{n-i}$. Hence, 
taking derivatives $n-i$ times and evaluating at $h=0$, we get $f^{(n-i)}(\lambda_{\varphi(i)})=0$.
In other words, $\res(f,f^{(n-i)})=0$.
\end{proof}

\begin{remark}
Another way to express $\epsilon_{n-i}(\lambda_1-\lambda_{k},\dots,\lambda_n-\lambda_{k})$
is by using the fact that it is equal to $f^{(i)}(\lambda_k)/i!$, 
\[
f(x)=\sum_{j\ge 0}^n (-1)^{n-j}\epsilon_{n-j}(\lambda_1,\dots,\lambda_n)x^j\Longrightarrow
\frac{f^{(i)}(\lambda_k)}{i!}=
\sum_{j\ge i}^n \binom{j}{i}(-1)^{n-j}\epsilon_{n-j}(\lambda_1,\dots,\lambda_n)\lambda_k^{j-i}.
\]
Then, for $1\le i\le n-1$ and $1\le k\le n$,
\[
\epsilon_{n-i}(\lambda_1-\lambda_{k},\dots,\lambda_n-\lambda_{k})=
\sum_{j\ge i}^n \binom{j}{i}(-1)^{n-j}\epsilon_{n-j}(\lambda_1,\dots,\lambda_n)\lambda_k^{j-i}.
\]
\end{remark}

\begin{remark}
Given that the polynomials 
\[
\epsilon_{1}^{1}\dots\epsilon_{1}^{n},\quad\dots,\quad
\epsilon_{n-1}^{1}\dots\epsilon_{n-1}^{n}
\]
are symmetric in $\{\lambda_1,\dots,\lambda_n\}$ we can rewrite them as polynomials in
the elementary symmetric polynomials $\{\epsilon_0,\dots,\epsilon_{n-1}\}$, or equivalently, in
the coefficients of $f=x^n+a_{n-1}x^{n-1}+\dots+a_0$. The resulting polynomials are precisely the resultants,
\[
\res(f,f^{(i)}/i!),\quad 1\le i\le n-1.
\]
Indeed, the determinant of $f^{(i)}(F)/i!$, where $F$ is the companion matrix of $f$, is equal to 
$\res(f,f^{(i)}/i!)$, but also, is equal to $f^{(i)}(\lambda_1)/i!\dots f^{(i)}(\lambda_n)/i!$.
\end{remark}

\section{Abel-Gontcharoff polynomials}

In this section we review some known results about Abel-Gontcharoff polynomials
that can be found in \cite{levinson}. 
\begin{definition}
The Abel-Gontcharoff polynomial of degree $0$ is $G(x) = 1$ and
the Abel-Gontcharoff polynomial of degree $n\ge 1$ is defined as 
\[
G(x;y_0,\dots,y_{n-1}) = (-1)^n n!\int_{y_0}^x dt_1\int_{y_1}^{t_1} dt_2\dots
\int_{y_{n-1}}^{t_{n-1}} dt_n,\quad n\ge 1.
\]
The normalizing factor of $(-1)^n n!$ is not present in \cite{levinson}. 
\end{definition}

\begin{lemma}\label{identities}.
Let $G(x;y_0,\dots,y_{n-1})\in\ZZ[x,y_0,\dots,y_{n-1},c]$ be the degree $n\ge 1$
Abel-Gontcharoff polynomial
and $c$ a new variable. Then,
\begin{itemize}
\item $G(y_0;y_0,\dots,y_{n-1})=0$.
\item $G(x;y_0,\dots,y_{n-1})$ is homogeneous of degree $n$
and 
\[
G(x+c;y_0+c,\dots,y_{n-1}+c)=G(x;y_0,\dots,y_{n-1}).
\]
\item 
\[
\frac{1}{k!}
\frac{\partial^k G(x;y_0,\dots,y_{n-1})}{\partial x^k} =
\binom{n}{k}
G(x;y_k,\dots,y_{n-1})
\]
In particular, if $f(x) = G(x;y_0,\dots,y_{n-1})$, then $f^{(k)}(y_k) = 0$.
\item $G(x;y_0,\dots,y_{n-1})$ is a monic polynomial in $x$, 
\item 
\[
\frac{\partial G(x;y_0,\dots,y_{n-1})}{\partial y_k} = -
n\binom{n-1}{k}G(x;y_0,\dots,y_{k-1})G(y_k;y_{k+1},\dots,y_{n-1})
\]
\end{itemize}
\end{lemma}
\begin{proof}
The first assertion follows from the definition. The
second, third and fifth assertion are proven in \cite[\S2,\S3]{levinson}. The fourth assertion
follows by computing $\partial^n G/\partial x^n$ and dividing by $n!$.
\end{proof}

\begin{corollary}\label{resta}
Let $G(x;y_0,\dots,y_{n-1})\in\ZZ[x,y_0,\dots,y_{n-1},c]$ be the degree $n\ge 1$
Abel-Gontcharoff polynomial, $c$ a new variable and let $0\le k<n$. Then,
\[
G(x;y_0,\dots,y_{n-1})-
G(x;y_0,\dots,y_{k-1},c,y_{k+1},\dots,y_{n-1})=
\binom{n}{k}
G(x;y_0,\dots,y_{k-1})G(c;y_k,\dots,y_{n-1}) .
\]
\end{corollary}
\begin{proof}
The result follows by integrating an identity from Lemma \ref{identities},
\begin{align*}
 G(x;y_0,\dots,y_{n-1})  
 &-
G(x;y_0,\dots,y_{k-1},c,y_{k+1},\dots,y_{n-1})\\
&=
\int_{c}^{y_k}
\frac{\partial G(x;y_0,\dots,y_{n-1})}{\partial y_k}  dy_k\\
&=
-\int_{c}^{y_k}
n\binom{n-1}{k}G(x;y_0,\dots,y_{k-1})G(y_k;y_{k+1},\dots,y_{n-1})dy_k\\
&=
-n\binom{n-1}{k}G(x;y_0,\dots,y_{k-1})
\int_{c}^{y_k}
G(y_k;y_{k+1},\dots,y_{n-1})dy_k \\
&=
-\frac{n}{n-k}\binom{n-1}{k}G(x;y_0,\dots,y_{k-1})
\left(
G(y_k;y_k,y_{k+1},\dots,y_{n-1})-
G(c;y_k,y_{k+1},\dots,y_{n-1})
\right)
\\
&=
\binom{n}{k}G(x;y_0,\dots,y_{k-1})
G(c;y_k,y_{k+1},\dots,y_{n-1}).
\end{align*}
\end{proof}

\begin{remark}
From Corollary \ref{resta} it follows that if $\lambda$ is a root of 
$f=G(x;y_0,\dots,y_{n-1})$, then we have the divisibility
\[
f^{(k)}(c) \,\big|\,
G(\lambda;y_0,\dots,y_{k-1},c,y_{k+1},\dots,y_{n-1})
\]
as polynomials in $c$.
\end{remark}

\begin{proposition}\label{igualdad}
Let $G(x;y_0,\dots,y_{n-1})\in\ZZ[x,y_0,\dots,y_{n-1},c_0,\dots,c_{n-1}]$ be the degree $n\ge 1$ Abel-Gontcharoff polynomial. Then,
\[
G(x;y_0,\dots,y_{n-1}) = G(x;c_0,\dots,c_{n-1}) + \sum_{i=0}^{n-1}\binom{n}{i}
G(x;c_0,\dots,c_{i-1})
G(c_i;y_i,\dots,y_{n-1}).
\]
Furthermore, 
\[
\frac{(n-k)!}{n!}
\frac{\partial^k G(x;y_0,\dots,y_{n-1})}{\partial x^k} = 
G(x;c_k,\dots,c_{n-1}) + \sum_{i=k}^{n-1}\binom{n-k}{i-k}
G(x;c_k,\dots,c_{i-1})
G(c_i;y_i,\dots,y_{n-1}).
\]
\end{proposition}
\begin{proof}
Applying Corollary \ref{resta} for $0\le i<n$ we get,
\begin{align*}
G(x;c_0,\dots,c_{i-1},y_i,\dots,y_{n-1}) &-
G(x;c_0,\dots,c_{i-1},c_i,y_{i+1}\dots,y_{n-1}) \\
&=
\binom{n}{i}
G(x;c_0,\dots,c_{i-1})G(c_i,y_i,\dots,y_{n-1}) .
\end{align*}
The first formula follows by taking an alternated sum. The formula for the derivative
is obtained by using the identities from Lemma \ref{identities} and the identity
\[
\binom{n}{k}\binom{n-k}{i-k} = \binom{i}{k}\binom{n}{i}.
\]
\end{proof}

\begin{remark}
Taking $c_0=\dots=c_{n-1}=0$ in the previous Proposition, we recover the standard
writing,
\[
G(x;y_0,\dots,y_{n-1}) = x^n + \sum_{i=0}^{n-1}\binom{n}{i}
x^{i}b_i, \quad b_i=G(0;y_i,\dots,y_{n-1}).
\]
\end{remark}

\section{The conjecture}

\begin{definition}
Let us define three different schemes over $\ZZ$.
The root schemes (\cite[\S 5.1]{frutos}), 
the coefficients scheme (\cite[\S 2]{infinity}) and the Abel-Gontcharoff  scheme.
The \emph{root scheme} $\mathcal{R}\subseteq\PP^{n-1}$
is defined by the ideal
\[
\langle 
\epsilon_{n-i}^1(\lambda_1,\dots,\lambda_n)
\dots\epsilon_{n-i}^n(\lambda_1,\dots,\lambda_n)
,\quad i=1,\dots,n-1
\rangle\subseteq\ZZ[\lambda_1,\dots,\lambda_{n}],
\]
the \emph{Abel-Gontcharoff  scheme} 
$\mathcal{G}\subseteq\PP^{n-1}$
is defined by the ideal 
\[
\langle
G(y_k;y_0,\dots,y_{n-1} ),\quad k=1,\dots,n-1
\rangle\subseteq\ZZ[y_0,\dots,y_{n-1}]
\]
and the \emph{coefficients scheme} $\mathcal{X}\subseteq\PP(n,n-1,\dots,1)$
is a weighted projective scheme defined by the ideal
\[
\langle
\res(P,H_k(P)),\quad k=1,\dots,n-1
\rangle\subseteq\ZZ[a_0,\dots,a_{n-1}]
\]
where $P=x^n+a_{n-1}x^{n-1}+\dots+a_1x+a_0\in\ZZ[a_0,\dots,a_{n-1}][x]$.

Notice that any point in
$\mathcal{R}(\CC)$ or in
$\mathcal{G}(\CC)$ or in
$\mathcal{X}(\CC)$ defines a polynomial $f$
such that $f^{(k)}( y_k )= f(y_k)=0$ for all $k=0,\dots,n-1$, that is, a CA-polynomial.
\end{definition}

\begin{proposition}
The dimensions of 
$\mathcal{R}$,
$\mathcal{G}$ and
$\mathcal{X}$ 
coincide.
\end{proposition}
\begin{proof}
Let us define two maps 
$\phi_1:\mathcal{R}\to\mathcal{X}$ 
and
$\phi_2:\mathcal{G}\to\mathcal{X}$,
\begin{align*}
\phi_1(\lambda_1:\dots:\lambda_n) :=&
\left(
\epsilon_{n}(\lambda_1,\dots,\lambda_n):\dots:\epsilon_1(\lambda_1,\dots,\lambda_n)
\right),
\\
\phi_2(y_0:\dots:y_{n-1}):=&
\left(
\binom{n}{0}G(0;y_0,\dots,y_{n-1}):\dots:
\binom{n}{n-2}G(0;y_{n-2},y_{n-1}):
\binom{n}{n-1}G(0;y_{n-1})
\right).
\end{align*}
The result follows by noting that these two maps are finite and surjective.
\end{proof}

\begin{proposition}
The following are equivalent
\begin{enumerate}
\item The Casas-Alvero conjecture is true.
\item The cardinal of $\mathcal{R}(\CC)$ is $1$.
\item The dimension of $\mathcal{R}(\CC)$ is zero.
\item $\{\epsilon_i^{1}\dots\epsilon_{i}^{n}\colon 1\le i\le n-1\}$
is a regular sequence over $\CC[\lambda_1,\dots,\lambda_n]$.
\item $\mathcal{R}(\CC)$ is finite.
\end{enumerate}
\end{proposition}
\begin{proof}
Let us prove $1$ implies $2$ and $5$ implies $1$. 
The other implications $i$ implies $i+1$ are clear, $2\le i\le 4$.
Assume the Casas-Alvero conjecture is true
and let $(\lambda_1:\dots:\lambda_n)\in \mathcal{R}(\CC)$. 
Then, from Proposition \ref{ca-equiv}, the polynomial
$(x-\lambda_1)\dots(x-\lambda_{n})$ is a CA-polynomial and 
from the hypothesis, $(\lambda_1:\dots:\lambda_n)=(1:\dots:1)$.

Assume now that $\mathcal{R}(\CC)$ is finite and let $f$
be a CA-polynomial of degree $n$. Then, 
from Proposition \ref{ca-equiv}
the roots of $f$ determine
a point $(\lambda_1:\dots:\lambda_n)\in\mathcal{R}(\CC)$.
But $f(x+b)$ is also a CA-polynomial for all $b\in\CC$.
Then, 
\[
(\lambda_1-b:\dots:\lambda_{n}-b)\in\mathcal{R}(\CC),\quad\forall b\in\CC.
\]
From the finiteness of $\mathcal{R}(\CC)$, it follows necessarily
$\lambda_1=\dots=\lambda_{n}$.
\end{proof}

\begin{remark}
From the previous proof we deduce that $\mathcal{R}(\CC)$ is a cone with vertex $(1:\dots:1)$. In other words, 
if $p\in\mathcal{R}(\CC)$, then the line joining $p$ and $(1:\dots:1)$ is contained in $\mathcal{R}(\CC)$.
Indeed, if $(s:t)\in\PP^1(\CC)$ and $t\neq 0$,
\[
(s+tp_0:\dots:s+tp_n)=(s/t+p_0:\dots:s/t+p_n)\in\mathcal{R}(\CC).
\]
The last containment follows from the fact that if $f=(x-p_0)\dots(x-p_n)$ is a CA-polynomial, 
then $f(x-s/t)$ is also a CA-polynomial. The same argument implies 
that $\mathcal{G}(\CC)$,  $\mathcal{G}(K)$ and $\mathcal{R}(K)$
are also cones for any field $K$.
\end{remark}


\section{First result}\label{first}
In this section we prove that if $n=p^{s+r}+p^s$ or 
$n=p^{s+r}+2p^s$ then $\dim\mathcal{R}(\CC)\le 1$.

\begin{lemma}\label{binoms}
Let $K$ be a field of characteristic $p$ and let $n=p^r+1$. Then,
\[
\binom{n}{k}=0,\quad 2\le k\le n-2.
\]
\end{lemma}
\begin{proof}
Notice that $n$ in base $p$ is equal to $n=1+0p+\dots+0p^{r-1}+1p^r$. Then,
\[
\binom{n}{k}=\binom{n_0}{k_0}\binom{n_1}{k_1}\dots\binom{n_r}{k_r}
=\binom{1}{k_0}\binom{0}{k_1}\dots\binom{0}{k_{r-1}}\binom{1}{k_r}.
\]
Given that $2\le k\le n-2$, it follows that $k_r=0$ and also that
$k_0>1$ or some $k_1,\dots,k_{r-1}>0$. In either case, the result follows.
\end{proof}

\begin{definition}
Following \cite{frutos}, 
let us define the $k$-th \emph{net derivative} of $f$
as $N_k(f):=H_k(f)/C(n,k)$, where $C(n,k)$ is the binomial coefficient
\[
C(n,k) :=\binom{n}{k}.
\]
From the equality $C(j,k)C(n,j)=C(n,k)C(n-k,j-k)$
and the relations $a_i=C(n,i)b_i$, $i=1,\dots,n-1$, 
it follows $N_k(f)\in\ZZ[b_k,\dots,b_{n-1}]$,
\begin{align*}
H_k(f)=&
\binom{n}{k}x^{n-k}+\binom{n-1}{k}a_{n-1}x^{n-1-k}+\dots+
\binom{k+1}{k}a_{k+1}x+a_k\\
=&
\binom{n}{k}x^{n-k}+\binom{n-1}{k}\binom{n}{n-1}b_{n-1}x^{n-1-k}+\dots+
\binom{k+1}{k}\binom{n}{k+1}b_{k+1}x+\binom{n}{k}b_k\\
=&
\binom{n}{k}\left(x^{n-k}+\binom{n-k}{n-k-1}b_{n-1}x^{n-1-k}+\dots+
\binom{n-k}{k+1-k}b_{k+1}x+b_k\right)\\
=&
\binom{n}{k}N_k(f).
\end{align*}
\end{definition}


\begin{definition}
Consider the field extension $\mathbb{Q}\subseteq \mathbb{C}$
and let us choose an absolute value $|\cdot|$ on $\mathbb{C}$ compatible with the $p$-adic absolute value in $\mathbb{Q}$.
Recall that the closed ball $R=\{x\in\mathbb{C}\,\colon\,|x|\le 1\}$ is a ring with maximal ideal
the open ball $\mathfrak{m}=\{x\in\mathbb{C}\,\colon\,|x|<1\}$. The field $K=R/\mathfrak{m}$ is
of characteristic $p$ and $K=\overline{K}$, \cite[XII]{lang}. 
We say that a monic CA-polynomial $f$ of degree $n$
is written in \emph{normal form} if 
all its roots are in $R$ and  if it is written as
\[
f=x^n-nx^{n-1}+\binom{n}{n-2}b_{n-2}x^{n-2}+\dots+\binom{n}{2}b_2x^2,
\]
where $|b_k|\le 1$ for $k=2,\dots,n-2$.
According to \cite[(3.4)]{castryck}, by applying the transformations $f(x-a)$ and $f(b x)/b^n$, we can convert any monic CA-polynomial to 
a polynomial in normal form.
\end{definition}


\begin{proposition}
If $\dim\mathcal{R}(K) = d$, then $\dim\mathcal{R}(\CC)\le d$. 
\end{proposition}
\begin{proof}
By hypothesis, the fiber of $\mathcal{R}$ at $p$ has dimension $d$. Given that
$\spec(\ZZ)$ has dimension $1$, it follows from  \cite[\S 4.3.1, Th.3.12]{liu}
that $\dim\mathcal{R}\le d+1$.
Hence, by generic flatness \cite[Th.6.9.1]{ega4.2}, the
dimension of the generic fiber is less than or equal to $d$,
$\dim\mathcal{R}(\QQbar)\le d$.
Finally, let us prove $\dim\mathcal{R}(\QQbar)=\dim\mathcal{R}(\CC)$
by using Noether normalization lemma 
\cite[\S 5 Ex.16]{atiyah}
and the fact that $\CC$ is flat over $\QQbar$.
If $A$ is a finitely generated $\QQbar$-algebra, then there exists
a polynomial ring $P$ such that $P\subseteq A$ and
$A$ is a finitely generated $P$-module. 
Tensoring by $\CC$, it follows
that $P_\CC\subseteq A_\CC$ and $A_\CC$ is a finitely generated $P_\CC$-module.
Hence, 
\[
\dim(A)=\dim(P)=\dim(P_\CC)=\dim(A_\CC).
\]
\end{proof}

\begin{lemma}\label{en_K}
Any CA-polynomials of degree $n=p^r+1$
in normal form has one simple root equal to $1$
and the other roots in $\mathfrak{m}$.
\end{lemma}
\begin{proof}
The class of $f=x^n-nx^{n-1}+C(n,n-2)b_{n-2}x^{n-2}+\dots+C(n,2)b_2x^2$
in $K$ is equal to $(x-1)x^{n-1}$.
\end{proof}

\begin{lemma}\label{en_K2}
Assume $p>2$.
The only CA-polynomial
in normal form of degree $n=p^r+2$ in $K[x]$ is 
\[
x^{n-2}(x-1)^2.
\]
\end{lemma}
\begin{proof}
As in Lemma \ref{binoms}, $C(n,k)=0$ for $3\le k\le n-3$, 
$C(n,2)=C(n,n-2)=1$ and $C(n,n-1)=2$.
Then, any CA-polynomial in normal form over $K[x]$ can be written as
\[
f=x^n - 2x^{n-1} + b_{n-2}x^{n-2} + 
b_2x^2. 
\]
Two of its net derivatives are
$N_2(f)=x^{n-2} + b_2$ and 
$N_{n-2}(f)= x^2-2x+b_{n-2}$.
In order to compute the resultants, let us factorize the
derivatives. If $\alpha^{n-2}=b_2$, then $N_2(f)=(x+\alpha)^{n-2}$
and it follows $\res(f,N_2(f))=f(-\alpha)^{n-2}$.
Given that $N_{n-1}(f)=x-1$, we have $\res(f,N_{n-1}(f))= f(1)$.
Finally, if $N_{n-2}(f)=(x-1-\beta)(x-1+\beta)$, then
$\res(f,N_{n-2}(f))= f(1+\beta)f(1-\beta)$.
Summing up,
\[
\res(f,N_{2}(f))= b_2^{n-2}b_{n-2}^{n-2} + 2b_2^{n-1},\quad
\res(f,N_{n-2}(f))= b_2^2b_{n-2}^2,\quad
\res(f,N_{n-1}(f))= -1+b_{n-2}+b_2.
\]
Then, the only possibility is $f=x^n-2x^{n-1}+x^{n-2}$.
\end{proof}

\begin{theorem}\label{fiber}
The following statements are true
\begin{enumerate}
\item If there exists a finite number of CA-polynomials 
in normal form of degree $n$ 
in $K[x]$, then $\dim\mathcal{R}(\CC)\le 1$. 
\item Let $r$ be a positive integer. If $n=p^r+1$, then $\dim\mathcal{R}(\CC)\le 1$.
\item Let $r,s$ be a positive integers. If $n=p^{s+r}+p^s$, 
then $\dim\mathcal{R}(\CC)\le 1$.
\item Let $r$ be a positive integer. If $n=p^r+2$, then $\dim\mathcal{R}(\CC)\le 1$.
\item Let $r,s$ be a positive integers. If $n=p^{s+r}+2p^s$, 
then $\dim\mathcal{R}(\CC)\le 1$.
\end{enumerate}
\end{theorem}
\begin{proof}
The hypothesis of the first statement 
implies $\dim\mathcal{R}(K)=1$. 
The second statement follows from Lemma \ref{en_K}
which implies that there is only one
CA-polynomial in $K[x]$. 
For the third statement 
we deduce from  \cite[Prop.26]{infinity} (or \cite[Prop.9]{draisma}), 
that there is also one CA-polynomial in $K$ equal to $((x-1)x^{p^r})^{p^s}$.
The last two statements follow from Lemma \ref{en_K2}.
\end{proof}

\begin{remark}
For $n=p^{s+r}+p^s$ or 
$n=p^{s+r}+2p^s$, we deduce 
from Theorem \ref{fiber} 
that
if $\mathcal{R}(\CC)$ is not a point, then it is a finite union of lines joined at $(1:\dots:1)$.
Hence, we can apply
this result to several of the cases listed in \cite[6.5]{castryck},
\[
\begin{matrix}
20, 24, 28, 30, 35, 36, 40, 42, 45, 48, 55, 56, 60, \\
63, 66, 70, 72, 77, 78, 80, 84, 88, 90, 91, 98, 99, 100.
\end{matrix}
\]
For example, if $n\in\{20,24,28,30,40,56,66,35,45,99\}$,
then $\dim\mathcal{R}(\CC) \le 1$. Indeed,
\[
\begin{matrix}
20 = 2^ 2 + 2^ 4 = 1+19,&
24 = 1+23=2^3+2^4,&
28 = 1+3^3,\\
30 = 3+3^3 = 5+5^2 = 1+29,&
40=2^3+2^5,&
56=7+7^2,\\
66=2+2^6,&
35 = 2 . 5 + 5^2,&
45 = 2 . 3^2 + 3^3 = 2 + 43,\\
99 = 2 + 97 = 2 .3^2 + 3^4.
\end{matrix}
\]
This result bounds the dimension 
of all the cases listed except for $n\in\{70, 77, 78, 88, 100\}$.
\end{remark}

\section{Second result}\label{second}

In this section we prove that 
the coefficients of a CA-polynomial of degree $n=p^r+1$
in normal form are in $\QQbar$.

\begin{lemma}
Let $f\in\mathbb{C}[x]$ be a monic CA-polynomial of degree $n$ 
with roots $\lambda_1,\dots,\lambda_n$
and 
let $\sigma$ be a complex automorphism over $\mathbb{Q}$, \cite{kestelman}.
Then, $f^\sigma$ is another monic CA-polynomial of the same degree
with roots $\sigma(\lambda_1),\dots,\sigma(\lambda_n)$. Clearly, if 
$\lambda\in\mathbb{Q}$, then $\sigma(\lambda)=\lambda$.
\end{lemma}
\begin{proof}
If $f=x^n+a_{n-1}x^{n-1}+\dots+a_1x+a_0$, then
$f^\sigma = x^n+\sigma(a_{n-1})x^{n-1}+\dots+\sigma(a_1)x+\sigma(a_0)$.
Then,
\[
H_k(f^\sigma)=H_k(f)^\sigma=
\binom{n}{k}x^{n-k}+\binom{n-1}{k}\sigma(a_{n-1})x^{n-1-k}+\dots+
\binom{k+1}{k}\sigma(a_{k+1})x+\sigma(a_k).
\]
Now, if $\lambda$ is a root of $f$, it is easy to check that $\sigma(\lambda)$ is a root of $f^\sigma$,
\[
f^\sigma(\sigma(\lambda)) =
\sum_{i=0}^n \sigma(a_i) \sigma(\lambda)^i=
\sigma\left(\sum_{i=0}^n a_i \lambda^i\right) =
\sigma( f(\lambda) ) = 0.
\]
Finally, if $f$ is a CA-polynomial in normal form, it follows 
from above that
$f^\sigma$ is also a CA-polynomial in normal form.
\end{proof}

\begin{theorem}
Let $f\in\CC[x]$ be a CA-polynomial of degree $n=p^r+1$ in normal form.
Then, $f\in\QQbar[x]$.
\end{theorem}
\begin{proof}
Let $\{\lambda_1,\dots,\lambda_n\}$ be the roots of $f$, 
where $\lambda_1=1$ and $\lambda_2,\dots,\lambda_n\in\mathfrak{m}$.
Let $\sigma$ be a complex automorphism and let $\lambda\in\mathfrak{m}$ be a root such that
$|\sigma(\lambda)|=1$. Then $f^\sigma$ 
is another CA-polynomial in normal form without linear term and with two roots of absolute value $1$. This contradicts Lemma \ref{en_K}. 
Then, if some root of $f$ is transcendental, then there exists a complex
automorphism sending it to a transcendental of absolute value $1$ which is not
possible.
\end{proof}

\section{Third result}\label{third}
In this final section, we prove that there are no
CA-polynomial of degree $20$ with three recycled roots.

\begin{lemma}\label{newton}
If $y_0,\dots,y_{k-1}\in R$, then $G(x;y_0,\dots,y_{k-1})\in R[x]$ 
and its roots are in $R$.
\end{lemma}
\begin{proof}
Given that the coefficients of $g=G(x;y_0,\dots,y_{k-1})$ are polynomials in 
$y_0,\dots,y_{k-1}$, it follows that $g\in R[x]$. Then, 
from the Newton polygon of $g$, it follows that all of its roots
are also in $R$.
\end{proof}

\begin{lemma}\label{bound}
Let $y_0,\dots,y_{k-1},c_0,\dots,c_{k-1}\in R$ be such that
$y_i-c_i\in\mathfrak{m}$ for all $i=0,\dots,k-1$.
Consider the set of indices $S=\{i\,\colon\, c_i\neq y_i\}\subseteq\{0,\dots,k-1\}$,
the value $M=\max\{ |C(k,i)|\,\colon\, i\in S\}$ and an element $\lambda\in R$.
Then,
\[
|G(\lambda;y_0,\dots,y_{k-1})| =M
\iff
|G(\lambda;c_0,\dots,c_{k-1})|=M.
\]
\end{lemma}
\begin{proof}
From Corollary \ref{igualdad} we have 
\[
G(x;y_0,\dots,y_{k-1})
 = G(x;c_0,\dots,c_{k-1}) + \sum_{i=0}^{k-1}\binom{k}{i}
G(x;c_0,\dots,c_{i-1})
G(c_i;y_i,\dots,y_{k-1}),
\]
For any $i\in\{0,\dots,k-1\}$, $G(c_i;y_i,\dots,y_{k-1})$ is congruent to
$G(y_i;y_i,\dots,y_{k-1})$ modulo $\mathfrak{m}$, but given that
$G(y_i;y_i,\dots,y_{k-1})=0$, it follows $|G(c_i;y_i,\dots,y_{k-1})|<1$.
Furthermore, if $i\not\in S$, we have $c_i=y_i$ and then $G(c_i;y_i,\dots,y_{k-1})=0$.
Hence, 
\[
G(x;y_0,\dots,y_{k-1}) = G(x;c_0,\dots,c_{k-1}) + \sum_{i\in S}\binom{k}{i}
G(x;c_0,\dots,c_{i-1})
G(0;y_{i},\dots,y_{k-1}).
\]
For $i\in S$, since $G(x,c_0,\dots,c_{i-1})\in R[x]$ and $\lambda\in R$, 
it follows $|G(\lambda,c_0,\dots,c_{i-1})|\le 1$.
Finally, the absolute value of the sum at $x=\lambda$ is less than $M$.
\end{proof}

The next  Corollary gives a similar result as in \cite[Th.2]{castryck}.
\begin{corollary}
Let $n=p^r+1$,
let $y_2,\dots.y_{n-2}\in \mathfrak{m}\cup\{1\}$ and let $f(x)=G(x;0,0,y_2,\dots,y_{n-2},1)$.
Let $c_i=0$ if $|y_i|<1$ and $c_i=1$ if $y_i=1$.
If $|G(1;c_0,\dots,c_{n-1})|=1/p$, then $f$ is not a CA-polynomial.
\end{corollary}
\begin{proof}
Lemma \ref{bound} implies $|f(1)|=1/p$, but if 
$f$ is a CA-polynomial, then $f(1)=0$, a contradiction.
\end{proof}

\begin{remark}
For the case $n=20$ and $p=19$, we checked all the $2^{17}$ cases for 
\[
|G(1;0,0,c_2,\dots,c_{n-2},1)|=\frac{1}{p}.
\]
We obtained 
that there are $6680$ possible CA-polynomials with two or more different roots.
Avoiding the cases $y_4=y_{16}=1$, $y_5=y_{10}=1$ and 
$y_{10}=y_{15}=1$ we reduced the list to $3125$
cases of possible CA-polynomials of degree $20$.
For each case, the conditions of being a CA-polynomial 
produce a system with one more equation than variables.
For example, the
most computationally intensive cases, have
$16$ equations with $15$ variables. These 4 possible CA-polynomials are
\begin{align*}
G(x&;0, 0,y_2,y_3, y_4, y_5, y_6, y_7, y_8, y_9, 1, y_{11}, 1, y_{13}, y_{14}, y_{15}, y_{16}, y_{17}, y_{18}, 1),\\
G(x&;0, 0,y_2,y_3, y_4, y_5, y_6, y_7, y_8, y_9, 1,1, y_{12}, y_{13}, y_{14}, y_{15}, y_{16}, y_{17}, y_{18}, 1),\\
G(x&;0, 0,y_2,y_3, y_4, 1, y_6, y_7, y_8, y_9, y_{10}, y_{11}, y_{12}, y_{13}, y_{14}, y_{15}, 1, y_{17}, y_{18}, 1),\\
G(x&;0, 0,y_2,y_3, 1, y_5, y_6, y_7, y_8, y_9, y_{10}, y_{11}, y_{12}, y_{13}, y_{14}, y_{15}, y_{16}, 1, y_{18}, 1).
\end{align*}
\end{remark}

The next Proposition gives a similar result as in \cite[Prop.20]{castryck}.
\begin{proposition}
Let $n=p^r+1$ and let $m$ be 
the minimum of $|G(1;0,0,c_2,\dots,c_{n-2},1)|$ for all $c_2,\dots,c_{n-2}\in\{0,1\}$,
\[
m =  \min\left\{
|G(1;0,0,c_2,\dots,c_{n-2},1)|
\,\colon\,c_2,\dots,c_{n-2}\in\{0,1\}
\right\}.
\]
Let $f=G(x;0,0,y_2,\dots,y_{n-2},1)$ be a polynomial with 
$|y_i|\le m$ or $y_i=1$ for $2\le i\le n-2$. Then $f$ is not
a CA-polynomial
\end{proposition}
\begin{proof}
First recall from \cite[Prop.6]{draisma} that $m>0$.
Let $c_i=0$ if $|y_i|<1$ and $c_i=1$ if $y_i=1$.
Notice that $G(x;y_i,\dots,y_{n-2},1)$ is a polynomial in $R[x]$ with roots in $R$ and
$y_i$ one of its roots. Hence, its constant term 
$b_i := G(0;y_i,\dots,y_{n-2},1)$ satisfy $|b_i|\le |y_i|$. 
Let us write $f$ as
\[
f(x) = G(x;c_0,\dots,c_{n-1}) + \sum_{i\in S}\binom{n}{i}
G(x;c_0,\dots,c_{i-1})b_i,
\]
where $S=\{i\,\colon\, c_i=0\}$, 
$|b_i|\le |y_i|\le m$ and $|C(n,i)|<1$ for all $i\in S$. 
Then, the absolute value of the sum is less than $m$ which implies
$|f(1)| = |G(1;c_0,\dots,c_{n-1})|\ge m$. 
\end{proof}

\begin{remark}
For the case $n=20$ and $p=19$, the number $m$ from the previous Proposition
can be computed relatively easy and it is equal to $5$.
Hence, in any CA-polynomial of degree $20$ there exists a common 
root of absolute value greater than $(1/p)^5$.
The same result holds for $n=24$ and $p=23$.
\end{remark}

The next Propositions generalizes \cite[Prop.16]{castryck}
\begin{proposition}
Let $f$ be a degree $n$ complex polynomial having two or more different roots.
Assume there exist
a prime $q$, a set $S$ and a root $\lambda$ such that $|C(n,i)|_q<1$ for all $i\in S$
and $f^{(i)}(\lambda)=0$ for all $i\not\in S$.
Then, $f$ is not a CA-polynomial.
\end{proposition}
\begin{proof}
By changing $x$ by $x-\lambda$ in $f$ and scaling it by the root of maximum 
absolute value (with respect to $q$) we can write $f$ as 
\[
f=x^{n}+\binom{n}{n-1}b_{n-1}x^{n-1}+\dots+
\binom{n}{1}b_{1}x,
\]
where $|b_i|_q\le1$ for all $i$ and by hypothesis $b_i=0$ for all
$i\not\in S$. Then, evaluating at the root $1$, we get 
\[
-1 = \sum_{i\in S}\binom{n}{i}b_{i}.
\]
By taking the absolute value on both sides of the equality, we arrive at a contradiction.
\end{proof}

\begin{remark}
For the case $n=20$, we can apply the previous Proposition
to deduce that 
there is no CA-polynomial $G(x;y_0,\dots,y_{19})$ such that
$y_{1}=y_{19}$ or $y_{4}=y_{16}$ or $y_5=y_{10}=y_{15}$.
\end{remark}

\begin{theorem}
There is no CA-polynomial of degree $20$ with a root of multiplicity $11$ or more
and there is no CA-polynomial of degree $24$ with a root of multiplicity $15$ or more.
\end{theorem}
\begin{proof}
Let $f$ be a CA-polynomial of degree $20$ with a root of multiplicity $11$ or more.
Then, $f=G(x;y_0,\dots,y_{19})$ where we can assume $y_0=\dots=y_{10}=0$.
Hence, $y_5=y_{10}=0$.
Analogously, if $f$ is a 
CA-polynomial of degree $24$ with a root of multiplicity $15$ or more, 
then $y_7=y_ {14}=0$.
Let us prove that these two cases are not possible.
First case $n=20$ over $\mathbb{F}_5$. The radical of the ideal 
\[
\langle
f(1),\res(f,N_5(f)),\res(f,N_{10}(f)),\res(f,N_{15}(f))
\rangle
\]
gives three solutions,
\[
(b_5,b_{10},b_{15})\in \left\{
(0,-1,0),
(1,0,-2),
(-2,0,1)
\right\}.
\]
If $y_5=y_{10}=0$, then $b_5=b_{10}=0$ which is not possible.
Second case $n=24$ over $\mathbb{F}_7$. The radical of the ideal
\[
\langle 
f(1),\res(f,N_7(f)),\res(f,N_{14}(f)),\res(f,N_{21}(f))
\rangle
\]
gives the solutions
\[
(b_7,b_{14},b_{21})\in \left\{
(3,3,0),(-3,-1,0),(1,0,1),(0,-3,1),(-1,0,3),(0,2,-2)
\right\}.
\]
If $y_7=y_{14}=0$, then $b_7=b_{14}=0$ which is not possible.
\end{proof}

\begin{theorem}
There are no CA-polynomial of degree $20$ with three recycled roots.
That is, if $f=G(x;0,0,y_2,\dots,y_{n-2},1)$ and $y_i\in\{0,1,y\}$ for all $i$, then $f$
is not a CA-polynomial.
\end{theorem}
\begin{proof}
Let $i$ be a positive integer less than $3^{17}$
written in base $3$, $i=d_2+d_33+d_23^2+\dots+d_{18}3^{17}$
where $d_2,\dots,d_{18}\in\{0,1,2\}$.
Let $f_i$ be the polynomial defined as $f_i=G(x;0,0,y_2,\dots,y_{18},1)$ where
\[
y_k=
\begin{cases}
1&\text{ if } d_k=1\\
0&\text{ if } d_k=0\\
y&\text{ if } d_k=2
\end{cases},\quad k=2,\dots,18.
\]
Then, $f_i$ is a CA-polynomial if and only if $f_i(1)=f_i(y)=0$. In other words, 
if $\res(f_i(1),f_i(y))=0$ as a polynomial in $y$.
We checked each of the $3^{17}$ possible cases in less than 48 hours in a 
personal computer (4GHz of CPU and 4GB of memory). 

In order to avoid the computation of the resultant, we analyzed the slopes
of the Newton polygons of $f(1)$ and $f(y)$ (for several primes), but we concluded that
the computation of the resultant is faster.
\end{proof}


\end{document}